\documentclass[prd,showpacs,twocolumn]{revtex4-1}

\usepackage{amsmath}
\usepackage{amsfonts}
\usepackage{graphicx}

\newcommand{\RR}{{\mathbb R}}
\newcommand{\CC}{{\mathbb C}}
\newcommand{\HH}{{\mathbb H}}
\newcommand{\KK}{{\mathbb K}}
\newcommand{\OO}{{\mathbb O}}

\newcommand{\isom}{\cong}

\newcommand{\sa}{\mathfrak{sa}}
\newcommand{\so}{\mathfrak{so}}
\newcommand{\su}{\mathfrak{su}}
\renewcommand{\sl}{\mathfrak{sl}}
\renewcommand{\sp}{\mathfrak{sp}}
\renewcommand{\aa}{\mathfrak{a}}
\newcommand{\cc}{\mathfrak{c}}
\newcommand{\dd}{\mathfrak{d}}
\newcommand{\ee}{\mathfrak{e}}
\newcommand{\ff}{\mathfrak{f}}
\renewcommand{\gg}{\mathfrak{g}}
\newcommand{\der}{\mathfrak{der}}

\newcommand{\SO}{\hbox{SO}}
\newcommand{\SU}{\hbox{SU}}
\newcommand{\SL}{\hbox{SL}}
\newcommand{\Sp}{\hbox{Sp}}
\newcommand{\EE}{E}
\newcommand{\FF}{F}
\newcommand{\UU}{\mathcal{U}}

\newcommand{\eA}{A}
\newcommand{\eG}{G}
\newcommand{\eS}{S}
\newcommand{\eD}{D}
\newcommand{\eE}{E}
\newcommand{\eF}{F}
\newcommand{\eX}{X}
\newcommand{\eY}{Y}
\newcommand{\eZ}{Z}

\newcommand{\mat}[1]{\mathbf{#1}}
\newcommand{\Pmat}{\mat{P}}

\newcommand{\HQ}{\mathrm{H}_Q}
\newcommand{\HHH}{\mathrm{H}_3(\OO)}
\newcommand{\AAA}{{\cal A}}
\newcommand{\BBB}{{\cal B}}
\newcommand{\III}{{\cal I}}
\newcommand{\PPP}{{\cal P}}

\renewcommand{\Re}{\mathrm{Re}\,}
\renewcommand{\Im}{\mathrm{Im}\,}
\renewcommand{\bar}[1]{\overline{#1}}

\begin{document}

\title{\boldmath A New Division Algebra Representation of $\EE_6$}

\author{Tevian Dray}
\email{tevian@math.oregonstate.edu}
\affiliation{Department of Mathematics, Oregon State University,
Corvallis, OR  97331, USA}

\author{Corinne A. Manogue}
\email{corinne@physics.oregonstate.edu}
\affiliation{Department of Physics, Oregon State University,
Corvallis, OR  97331, USA}

\author{Robert A. Wilson}
\email{r.a.wilson@qmul.ac.uk}
\affiliation{School of Mathematical Sciences, Queen Mary,
University of London, London E1 4NS, UK}

\begin{abstract}
We construct the well-known decomposition of the Lie algebra $\ee_8$ into
representations of $\ee_6\oplus\su(3)$ using explicit matrix representations
over pairs of division algebras.  The minimal representation of $\ee_6$,
namely the Albert algebra, is thus realized explicitly within $\ee_8$, with
the action given by the matrix commutator in $\ee_8$, and with a natural
parameterization using division algebras.  Each resulting copy of the Albert
algebra consists of \textit{anti-Hermitian} matrices in $\ee_8$, labeled by
imaginary (split) octonions.  Our formalism naturally extends from the Lie
algebra to the Lie group $\EE_6\subset\EE_8$.
\end{abstract}


\maketitle

\section{Introduction}
\label{intro}

The exceptional Lie groups have long been suspected of playing a key role in
the modern description of fundamental
particles~\cite{Gunaydin,GurseyE6,GT,Bars80}, and the role of the octonions in
describing the exceptional Lie groups is well known.
In the early 1980s, there was an explosion of interest in using octonionic
structures to describe supersymmetry, starting
with~\cite{KugoTownsend83,OliveWest83}, and superstring
theory~\cite{FairlieI,Goddard,FairlieII,Corrigan,SudberyCAM,GSW}; a short
history is given in~\cite{BHsuperI}.
More recently, several
authors~\cite{Dixon,Furey,Furey18,Furey19,Furey22a,Furey22b,%
TDV,Krasnov21,Perelman21} have used division and/or Clifford algebra
structures to represent particle states.

Motivated by the considerable
recent interest in using the exceptional Lie group $\EE_8$ to describe
fundamental particles~\cite{Lisi,Chester,Magic,Octions}, we describe here a
mechanism for for analyzing both $\EE_6$ and its well-known action on the
Albert algebra in the language of the exceptional Lie algebra $\ee_8$.

Building on previous work by Sudbery~\cite{Sudbery}, Manogue and
Schray~\cite{Lorentz} gave an explicit construction of the Lie group
$\SU(2,\OO)$ as the double cover of $\SO(9)$, and of $\SL(2,\OO)$ as the
double cover of $\SO(9,1)$, and Manogue and Dray~\cite{Denver,York} (see
also~\cite{AaronThesis,Sub,Structure}) extended this construction to an
interpretation of $\SU(3,\OO)$ as $\FF_4$, and $\SL(3,\OO)$ as $\EE_{6(-26)}$.
More recently, Dray, Manogue, and Wilson~\cite{Cubies} further extended this
construction to $\EE_{7(-25)}$, providing an interpretation in terms of
$\Sp(6,\OO)$, thus providing an explicit realization at the Lie group level of
the description originally given at the Lie algebra level by Barton and
Sudbery~\cite{SudberyBarton}.  Each of the constructions above involves a
single copy of $\OO$, supporting the interpretation of the elements of the
first three rows of the Freudenthal--Tits magic square~\cite{Freudenthal,Tits}
as unitary, Lorentz, and symplectic groups over the division algebras $\RR$,
$\CC$, $\HH$, and $\OO$, respectively, as shown in Table~\ref{3x3}.

\begin{table*}
\small
\begin{center}
\begin{tabular}{|c|c|c|c|c|}
\hline
&$\RR$&$\CC$&$\HH$&$\OO$\\\hline
$\RR'$
&$\su(3,\RR)$&$\su(3,\CC)$&$\cc_3\isom\su(3,\HH)$&$\ff_4\isom\su(3,\OO)$\\
\hline
$\CC'$&
$\sl(3,\RR)$&$\sl(3,\CC)$&$\aa_{5(-7)}\isom\sl(3,\HH)$
  &$\ee_{6(-26)}\isom\sl(3,\OO)$\\
\hline
$\HH'$&
$\cc_{3(3)}\isom\sp(6,\RR)$&$\su(3,3,\CC)\isom\sp(6,\CC)$
  &$\dd_{6(-6)}\isom\sp(6,\HH)$&$\ee_{7(-25)}\isom\sp(6,\OO)$\\
\hline
$\OO'$&
$\ff_{4(4)}$&$\ee_{6(2)}$&$\ee_{7(-5)}$&$\ee_{8(-24)}$\\
\hline
\end{tabular}
\end{center}
\caption{The ``half-split'' real form of the Freudenthal--Tits magic square of
Lie algebras.}
\label{3x3}
\end{table*}

The identifications in Table~\ref{3x3} do not appear to admit a natural
extension to the fourth row of the magic square.  This result is not
surprising, since unitary, Lorentz, and symplectic groups are normally defined
in terms of the vectors on which they act, whereas the fourth row of the magic
square contains (the Lie algebra) $\ee_8$, whose minimal representation is the
adjoint representation.  In other words, there is no smaller space of
``vectors'' on which the ``matrices'' in $\ee_8$ can act; these matrices must
act on each other.  In order to provide a description of the fourth row, it is
therefore necessary to describe the adjoint representation itself in terms of
division algebras.
In~\cite{Magic,Octions}, we accomplished just that, providing a
description of the Lie algebra $\ee_8$---and hence the entire magic
square---in terms of $3\times3$ matrices over the tensor product of
\textit{two} division algebras.  This construction, which is summarized in
Section~\ref{e8}, provides an explicit implementation of the Vinberg
description~\cite{Vinberg} of the magic square.  A simpler variation of this
approach was previously used by Dray, Huerta, and
Kincaid~\cite{JoshuaThesis,SO42,2x2} to describe the analogous ``$2\times2$''
magic square, involving real forms of subgroups of $\SO(16)$, in terms of two
division algebras.

Working backward, it is straightforward to reduce the description of the Lie
algebra $\ee_8$ given in~\cite{Magic,Octions} to a description of the adjoint
representation of each of the other Lie algebras in the magic square.  But we
can do more: Breaking the adjoint representation of one algebra down to the
adjoint representation of a subalgebra results in a decomposition of the
original algebra into several pieces, and it is instructive to examine these
pieces explicitly.

In this paper, we analyze the well-known decomposition of $\ee_8$ into
representations of $\ee_6\oplus\su(3)$ from this perspective.  Although the
decomposition itself is well known, the resulting use of $\ee_8$ to express
the action of $\ee_6$ on the Albert algebra in terms of matrix commutators is,
to our knowledge, new.  We include a discussion of different possible real
forms, corresponding to the use of either split or regular division algebras.
In a separate paper~\cite{MagicE7}, we will provide a similar discussion of
the decomposition of $\ee_8$ into representations of $\ee_7\oplus\su(2)$.

\section{Background}
\label{background}

\begin{figure}
\centering
\includegraphics[width=5cm]{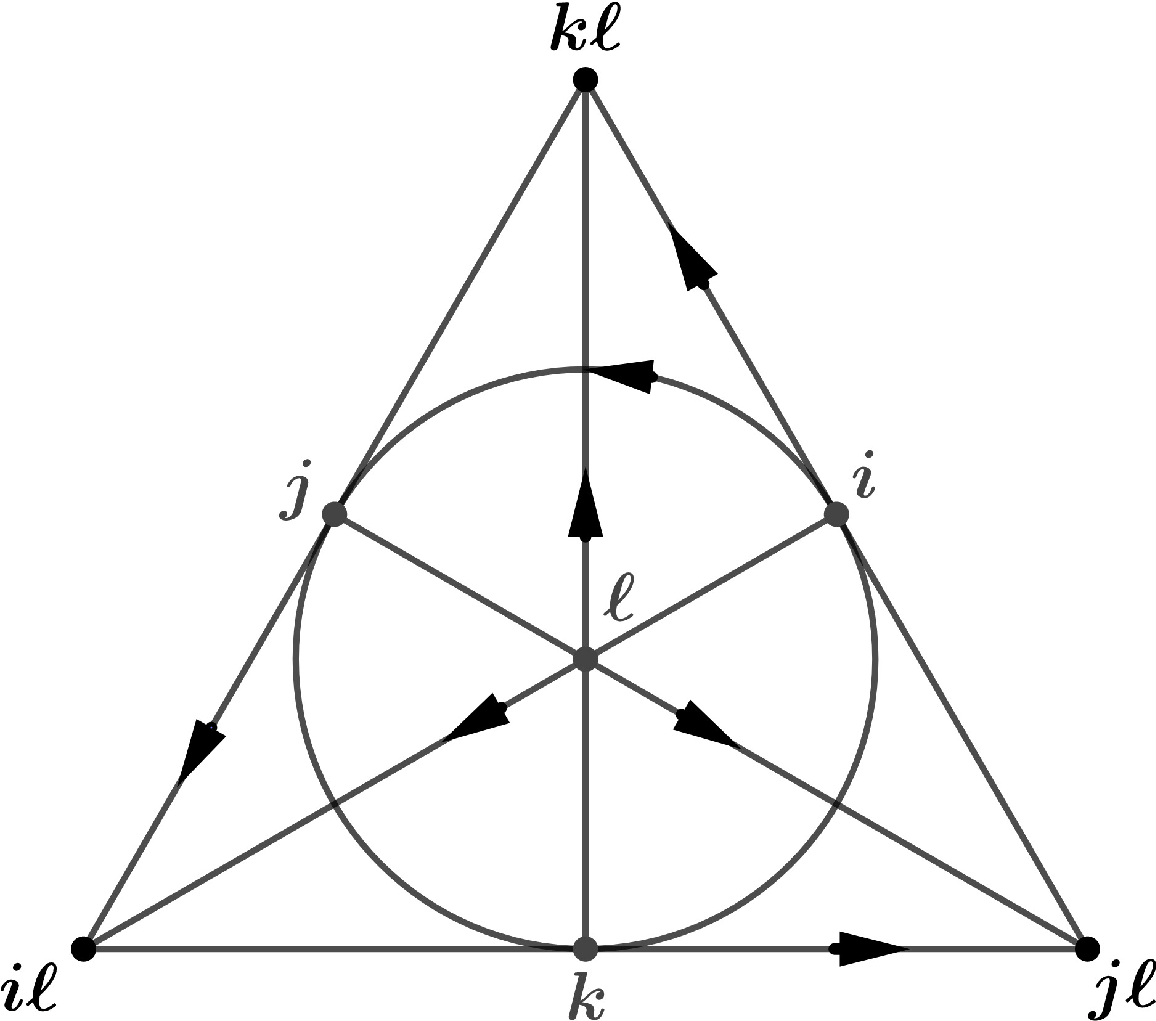}
\caption{A graphical representation of the  octonionic multiplication table.}
\label{omult3}
\end{figure}

\begin{table}
\centering
\small
\begin{tabular}[b]{|c|c|c|c|c|c|c|c|}
\hline
&\boldmath$i$&\boldmath$j$&\boldmath$k$&\boldmath$k\ell$
  &\boldmath$j\ell$&\boldmath$i\ell$&\boldmath$\ell$\\\hline
\boldmath$i$&$-1$&$k$&$-j$&$j\ell$&$-k\ell$&$\ell$&$i\ell$\\\hline
\boldmath$j$&$-k$&$-1$&$i$&$-i\ell$&$\ell$&$k\ell$&$j\ell$\\\hline
\boldmath$k$&$j$&$-i$&$-1$&$-\ell$&$i\ell$&$-j\ell$&$k\ell$\\\hline
\boldmath$k\ell$&$-j\ell$&$i\ell$&$\ell$&$-1$&$i$&$-j$&$-k$\\\hline
\boldmath$j\ell$&$k\ell$&$\ell$&$-i\ell$&$-i$&$-1$&$k$&$-j$\\\hline
\boldmath$i\ell$&$\ell$&$-k\ell$&$j\ell$&$j$&$-k$&$-1$&$-i$\\\hline
\boldmath$\ell$&$-i\ell$&$-j\ell$&$-k\ell$&$k$&$j$&$i$&$-1$\\\hline
\noalign{\vspace{0.15in}}
\end{tabular}
\caption{The octonionic multiplication table.}
\label{omult}
\end{table}

The \textit{octonions} $\OO$ are the real algebra spanned by the identity
element~$1$ and seven square roots of~$-1$ that we denote $i$, $j$, $k$,
$k\ell$, $j\ell$, $i\ell$, $\ell$ and refer to collectively as $\UU$, whose
multiplication table is neatly described by the oriented Fano geometry shown
in Figure~\ref{omult3}, and given explicitly in Table~\ref{omult}.  The
\textit{split octonions} $\OO'$ are the real algebra spanned by the identity
element~$1$, three square roots of~$-1$ that we denote $I$, $J$, $K$, and four
square roots of~$+1$ that we denote $KL$, $JL$, $IL$, $L$, whose
multiplication table is given in Table~\ref{smult}.  The multiplication table
for the split octonions can also be expressed using a modified Fano geometry
along the lines of Figure~\ref{omult3}, so long as one remembers to change the
sign of all products involving two imaginary units containing an $L$, as can
be seen by comparing the four $4\times4$ blocks in Tables~\ref{omult}
and~\ref{smult}.

\begin{table}
\centering
\small
\begin{tabular}{|c|c|c|c|c|c|c|c|}
\hline
$$&\boldmath$I$&\boldmath$J$&\boldmath$K$
  &\boldmath$KL$&\boldmath$JL$&\boldmath$IL$&\boldmath$L$\\\hline
\boldmath$I$&$-1$&$K$&$-J$&$JL$&$-KL$&$-L$&$IL$\\\hline
\boldmath$J$&$-K$&$-1$&$I$&$-IL$&$-L$&$KL$&$JL$\\\hline
\boldmath$K$&$J$&$-I$&$-1$&$-L$&$IL$&$-JL$&$KL$\\\hline
\boldmath$KL$&$-JL$&$IL$&$L$&$1$&$-I$&$J$&$K$\\\hline
\boldmath$JL$&$KL$&$L$&$-IL$&$I$&$1$&$-K$&$J$\\\hline
\boldmath$IL$&$L$&$-KL$&$JL$&$-J$&$K$&$1$&$I$\\\hline
\boldmath$L$&$-IL$&$-JL$&$-KL$&$-K$&$-J$&$-I$&$1$\\\hline
\end{tabular}
\caption{The split octonionic multiplication table.}
\label{smult}
\end{table}

Vinberg~\cite{Vinberg} gave a symmetric description of the Lie algebras in the
Freudenthal--Tits magic square.  In particular, the Vinberg description of
(half-split) $\ee_8$ is
\begin{align}
\ee_{8(-24)}
  &= \sa(3,\OO'\otimes\OO)\oplus\der(\OO')\oplus\der(\OO) 
\end{align}
where $\sa(\KK)$ denotes $3\times3$ tracefree, anti-Hermitian matrices over
$\KK$, and where $\der(\KK)$ denotes the derivations (infinitesimal
automorphisms) of $\KK$, so $\der(\OO)=\gg_2$, and $\der(\OO')=\gg_{2(2)}$.
There are $3\times8\times8=192$ off-diagonal elements and $2\times(7+7)=28$
diagonal elements in $\sa(3,\OO'\otimes\OO)$, and $14+14=28$ elements in
$\gg_2\oplus\gg_{2(2)}$, for a total of $192+28+28=248$ elements in $\ee_8$,
as expected.  The other Lie algebras in the Freudenthal--Tits magic square are
obtained by replacing $\OO$ and $\OO'$ by subalgebras; this labeling is shown
in the row and column headings in Table~\ref{3x3}.

\section{\boldmath The adjoint representation of $\ee_8$}
\label{e8}

We present here a summary of the construction of the adjoint representation of
$\ee_8$ given in~\cite{Magic}, yielding a description of $\ee_{8(-24)}$ as
\textit{generalized} $3\times3$ matrices over $\OO'\otimes\OO$, which we write
as $\su(3,\OO'\otimes\OO)$.  The construction proceeds in three steps.  In the
first two steps, we disentangle the action of certain subalgebras ($\so(9)$
and $\ff_4$) of $\ee_8$ on known (but non-adjoint) representations from the
representation, showing in each case that our generalized matrices themselves
correctly reproduce the adjoint action.  In the last step, we extend these
adjoint representations directly to $\ee_8$.

\paragraph{Step 1: $\so(9)$.}

To construct the adjoint representation of $\so(9)$, we first mimic the
construction of $\SO(9,1)$ in~\cite{Lorentz} by representing elements of the
vector representation $\RR^9$ of $\so(9)$ as $2\times2$ tracefree Hermitian
matrices $\Pmat$, so that
\begin{equation}
\Pmat = \begin{pmatrix}z& a\\ \bar{a}& -z\\\end{pmatrix}
\label{Pmat}
\end{equation}
with $a\in\OO$.  Now consider the off-diagonal anti-Hermitian matrices
\begin{equation}
\eX_q = \begin{pmatrix}0& q\\ -\bar{q}& 0\\\end{pmatrix}
\end{equation}
with $q\in\OO$, where $\eX_q$ acts on $\Pmat$ via
\begin{equation}
\Pmat \longmapsto \eX_q\Pmat + \Pmat \eX_q^\dagger = [\eX_q,\Pmat] .
\end{equation}
If we further impose $|q|=1$, then direct calculation shows that $\eX_q$ is
the infinitesimal generator of a (counterclockwise) rotation in the
$(q,z)$-plane, in agreement with~\cite{Lorentz}.

Given two operators $A$, $B$ acting on $\Pmat$, their commutator $[A,B]$ is
the operator given by
\begin{equation}
\Pmat \longmapsto \bigl[A,[B,\Pmat]\bigr] - \bigl[B,[A,\Pmat]\bigr] .
\label{commute}
\end{equation}
We can therefore construct $\so(9)$ by taking commutators of the
form~(\ref{commute}), and we define
\begin{align}
\eD_q &= \frac12 \> [\eX_1,\eX_q] ,\label{Dq}\\
\eD_{p,q} &= \frac12 \> [\eX_p,\eX_q] ,\label{Dpq}
\end{align}
where now $p,q\in\Im\OO$, and we further assume that $p\perp q$, that is
(using the polarized norm), that
\begin{equation}
p\bar{q}+q\bar{p}=0 ,
\end{equation}
and that $|p|=1=|q|$.  We will usually restrict $p$ and $q$ to be distinct
elements of our imaginary basis $\UU$.  Writing $x$ for $\Re a$
in~(\ref{Pmat}), it is easily seen (referring to~\cite{Lorentz} if necessary)
that $\eD_q$ generates a (counterclockwise) rotation in the $(x,q)$-plane, and
that $\eD_{p,q}$ generates a (counterclockwise) rotation in the $(p,q)$-plane.

The seven independent elements $\eD_q$ with $q\in\UU$, together with the 21
independent elements~$\eD_{p,q}$, with $p\ne q\in\UU$, form a basis of
$\so(8)$; the eight elements $\eX_q$ (where now $q\in\OO$) extend $\so(8)$ to
$\so(9)$.  A central idea in~\cite{Lorentz} is that the \textit{transverse
rotations} $\eD_{p,q}$ are \textit{nested}; they can not be implemented in
terms of a single matrix operation.

It is straightforward to verify that although $\eD_q$ is defined
by~(\ref{commute}) and~(\ref{Dq}) as a nested operation, it can nonetheless be
expressed using direct matrix multiplication as
\begin{equation}
\eD_q = \begin{pmatrix}q& 0\\ 0& -q\\\end{pmatrix} ;
\end{equation}
there are no associativity issues here.  The key to the matrix construction
in~\cite{Magic} is that the nested actions $\eD_{p,q}$ can also be represented
as \textit{generalized} $2\times2$ matrices, allowing us to dispense with the
vectors $\Pmat$ altogether, and compute commutators using matrices.  An
example of such a generalized matrix is
\begin{equation}
\eD_{p,q} = \begin{pmatrix}p\circ q& 0\\ 0& p\circ q\\\end{pmatrix}
\label{comp}
\end{equation}
where $\circ$ denotes \textit{composition}, that is
\begin{equation}
(p\circ q)r = p(qr)
\end{equation}
for $p,q,r\in\OO$.  Composition is associative by definition!  Furthermore,
composition respects anticommutativity, in the sense that
\begin{equation}
qp = -pq \Longrightarrow q\circ p=-p\circ q
\end{equation}
which holds in particular for $p\ne q\in\UU$.  For such $p$, $q$,
straightforward computation also shows that
\begin{align}
p\circ(p\circ q) &= p^2 q = -(p\circ q)\circ p \label{ppq} ,\\
(p\circ q)\circ(q\circ r) &= p\circ q^2 \circ r = q^2 p\circ r \label{pqr}
\end{align}
for any $r\in\OO$, since $p^2,q^2\in\RR$.

Computing commutators using matrix multiplication and $\circ$, we obtain the
same commutators as when using~(\ref{commute}), namely
\begin{align}
[\eX_1,\eX_p] &= 2\eD_p ,\label{x1xp}\\
[\eX_p,\eX_q] &= 2\eD_{p,q} ,\label{xpxq}\\
[\eX_1,\eD_p] &= -2\eX_p ,\label{x1dp}\\
[\eX_p,\eD_p] &= 2\eX_1 ,\label{xpdp}\\
[\eD_p,\eD_q] &= 2\eD_{p,q} ,\label{dpdq}\\
[\eX_p,\eD_{p,q}] &= -2\eX_q ,\label{xpdpq}\\
[\eD_p,\eD_{p,q}] &= -2\eD_q ,\label{dpdpq}\\
[\eD_{q,p},\eD_{p,r}] &= -2 \eD_{q,r} ,\label{dpqdqr}
\end{align}
where $p,q,r$ are distinct elements of $\UU$; all other commutators are zero.
We have therefore verified that the (generalized) matrix action of $\so(9)$
on $\RR^9$ is a Lie algebra homomorphism, that is, that the (generalized)
matrices themselves form a Lie algebra, namely the adjoint representation of
$\so(9)$.

This construction works just as well over $\OO'$ as over~$\OO$;
\eqref{x1xp}--\eqref{dpqdqr} still hold, although there are sign changes
in~\eqref{xpdp} and~\eqref{xpdpq}--\eqref{dpqdqr} due to an implicit factor of
$-p^2$, which is $1$ over $\OO$ but not over $\OO'$.

\paragraph{Step 2: $\ff_4$.}

Using a similar technique, to construct the adjoint representation of $\ff_4$,
we mimic the construction of $\ee_6$ in~\cite{Denver,York} (see
also~\cite{AaronThesis,Structure}) by considering the action of $\FF_4$ on the
Albert algebra $\HHH$ of $3\times3$ octonionic Hermitian matrices over $\OO$.
As in the $2\times2$ construction of $\so(9)$ above, the Lie algebra $\ff_4$
is generated by $3\times3$ octonionic anti-Hermitian matrices $\AAA$, acting
via
\begin{equation}
\PPP \longmapsto [\AAA,\PPP]
\end{equation}
on elements $\PPP\in\HHH$, where the commutator is computed using ordinary
matrix multiplication.  But $\FF_4$ is essentially the union of three
copies of $\SO(9)$.  We therefore begin by introducing the elements
\begin{align}
\eX_p = \begin{pmatrix}0& p& 0\\ -\bar{p}& 0 &0\\ 0& 0& 0\\\end{pmatrix} ,
\qquad
\eD_q = \begin{pmatrix}q& 0& 0\\ 0& -q& 0\\ 0& 0& 0\\\end{pmatrix} ,
\end{align}
with $p\in\OO$ and $q\in\Im\OO$, where we have recycled the names from our
description of $\so(9)$ in an obvious way; we refer to these elements as being
of \textit{type~I}.  There are corresponding elements with other block
structures, namely
\begin{align}
\eY_p &= \begin{pmatrix}0& 0& 0\\ 0& 0& p\\ 0& -\bar{p} &0\\\end{pmatrix} ,
\qquad
\eE_q = \begin{pmatrix}0& 0& 0\\ 0& q& 0\\ 0& 0& -q\\\end{pmatrix} ,
\end{align}
and
\begin{align}
\eZ_p &= \begin{pmatrix}0& 0& -\bar{p}\\  0& 0& 0\\ p& 0 &0\\\end{pmatrix} ,
\qquad
\eF_q = \begin{pmatrix}-q& 0& 0\\ 0& 0& 0\\ 0& 0& q\\\end{pmatrix} ,
\end{align}
which we refer to as \textit{type~II} and \textit{type~III}, respectively, as
well as their commutators
\begin{align}
\eD_{p,q} &= \frac12 \> [\eD_p,\eD_q] = \frac12 \> [\eX_p,\eX_q] ,\\
\eE_{p,q} &= \frac12 \> [\eE_p,\eE_q] = \frac12 \> [\eY_p,\eY_q] ,\\
\eF_{p,q} &= \frac12 \> [\eF_p,\eF_q] = \frac12 \> [\eZ_p,\eZ_q] .
\end{align}
As before, we normally assume that $p,q$ are distinct elements of $\UU$.

The 28 independent $\eD$s, $\eE$s, and $\eF$s each generate a copy of
$\so(8)$.  However, \textit{triality} tells us that these three copies of
$\so(8)$ must be identical!  Working out the action of these generators on
$\HHH$, we can indeed write all of the $\eE$s and $\eF$s in terms of the
$\eD$s, etc.; explicit formulas can be found in~\cite{AaronThesis,Magic} and
elsewhere (see also Appendix~\ref{so8}), but will not be needed here.

We therefore have 28 independent elements in $\so(8)$---the $\eD$s,
say---together with eight $\eX$s, eight $\eY$s, and eight $\eZ$s, for a total
of $28+3\times8=52$ elements, the correct number to generate $\ff_4$.  We now
show that these 52 elements do exactly that, thus verifying that the
commutators are in fact those of $\ff_4$.

We already know that the $\eD$s and $\eX$s together generate a copy of
$\so(9)$, using matrix commutators and $\circ$; similarly, so do the $\eE$s
and $\eY$s, or the $\eF$s and $\eZ$s.  However, matrix commutators between,
say, (double-index) $\eD$s and $\eE$s do \textit{not} always yield the correct
answer; the matrix operations do not respect triality.  We turn this bug into
a feature by requiring commutators between $\so(8)$ and other generators to be
computed \textit{using matrices of a single type}.  That is, \textit{before}
the commutator can be computed using (generalized) matrix multiplication,
triality must be used to express the element of $\so(8)$ as a matrix of the
correct type.  The resulting commutators are merely $3\times3$ versions of the
$\so(9)$ commutators in Step 1, and thus correctly reproduce the action on
$\HHH$.

It only remains to compute commutators of the $\eX$s, $\eY$s, and $\eZ$s with
each other.  But this is straightforward: Since the Albert algebra is a
representation of $\ff_4$, the Jacobi identity holds, which allows us to
compute
$[\eX_i,\eY_j]$ as the operator
\begin{align}
\AAA \longmapsto 
  &= \bigl[\eX_i,[\eY_j,\AAA]\bigr] - \bigl[\eY_j,[\eX_i,\AAA]\bigr]
\nonumber\\
  &= \bigl[[\eX_i,\eY_j],\AAA\bigr]
   = [\eZ_k,\AAA]
\end{align}
since the matrix commutator of $\eX_i$ and $\eY_j$ is $\eZ_k$.  We conclude
that we can use matrix commutators without impediment to compute Lie algebra
commutators of off-diagonal elements of~$\ff_4$.  Direct computation then
immediately yields
\begin{align}
[\eX_p,\eY_q] &= \eZ_{-\bar{pq}} ,\label{XYZ1}\\
[\eY_p,\eZ_q] &= \eX_{-\bar{pq}} ,\label{XYZ2}\\
[\eZ_p,\eX_q] &= \eY_{-\bar{pq}} ,\label{XYZ3}
\end{align}
for $p,q\in\OO$.

We have therefore shown that (generalized) matrix commutators---using $\circ$
when necessary, and converting type structure when appropriate---correctly
reproduce the Lie algebra structure of~$\ff_4$; it is no longer necessary to
act on the Albert algebra in order to compute commutators.  Further details
can be found in~\cite{Magic}.

\paragraph{Step 3: $\ee_8$.}

The hard work is done.  The minimal representation of $\ee_8$ is the adjoint
representation, so our matrix representation must now stand alone.  We need to
work over $\OO'\otimes\OO$, but everything we have done so far works just as
well over $\OO'$ as over $\OO$, since \eqref{XYZ1}--\eqref{XYZ3} are unchanged
over $\OO'\otimes\OO$, and the generalization of \eqref{x1xp}--\eqref{dpqdqr}
is straightforward.

To finalize the construction, consider first the diagonal elements.  An
anti-Hermitian element of $\OO'\otimes\OO$ must (be a sum of terms which) have
exactly one real and one imaginary factor, and hence lie in $\Im\OO'$
\textit{or} $\Im\OO$.  Since triality relates matrices of different type
structures, we only need to consider one type structure, and therefore have
$28+28=56$ independent diagonal elements, representing $\so(4,4)\oplus\so(8)$.

What's left?  There are now $3\times8\times8=192$ off-diagonal elements, whose
commutators can be computed using matrix multiplication---and $\circ$ where
necessary.  The only remaining subtlety is to check that the algebra closes,
but the anticommutativity of the imaginary units in each of $\OO'$ and $\OO$
ensures that the only diagonal elements we can obtain are those we already
have.

We're done; $28+28+192=248$, although our basis is slightly different than in
the counting argument given at the end of Section~\ref{background}.

\paragraph{Preferred basis.}

We have constructed three different bases of $\so(8)$, namely
$\{\eD_p,\eD_{p,q}\}$, $\{\eE_p,\eE_{p,q}\}$, $\{\eF_p,\eF_{p,q}\}$, of
types~I, II, and~III, respectively, which are related by triality.  However,
none of these bases are adapted to the preferred $\gg_2$ subalgebra of
$\so(8)$ that preserves the octonionic multiplication table.  Such a basis is
constructed in Appendix~\ref{so8}, with $\{\eA_p,\eG_p\}$ generating $\gg_2$,
$\{\eS_p\}$ being the missing generators of the $\so(7)$ that fixes $\Im\OO$,
and $\{\eD_p\}$ being the remaining generators of $\so(8)$, where $p\in\UU$.

\section{\boldmath Representations of $\ee_6$}
\label{e6}

\subsection{Overview}

It is well known that (complex) $\ee_8$ can be decomposed as
\begin{equation}
\ee_8 = \ee_6 \oplus \su(3) \oplus 162,
\end{equation}
since the centralizer (maximal commuting subalgebra) of $\ee_6$ in $\ee_8$ is
$\su(3)$.  The remaining $248-78-8=162$ elements of $\ee_8$ can be divided
into simultaneous representations of $\ee_6$ and $\su(3)$, normally written as
\begin{equation}
162 = 3\times27 + \bar{3}\times\bar{27} ,
\end{equation}
where each ``$3$'' denotes a minimal representation of $\su(3)$, each ``$27$''
denotes a minimal representation of $\ee_6$, and the bars denote the dual
representations.  Since the minimal representation of $\ee_6$ is the Albert
algebra, it must be the case that $\ee_8$ contains six copies of the Albert
algebra.  This decomposition holds for any real form of $\ee_8$, with minor
variations.  In the remainder of this section, we construct this decomposition
explicitly for $\ee_{6(-26)}\subset\ee_{8(-24)}$.

\subsection{\boldmath The structure of $\ee_6$}
\label{e6gen}

We can restrict $\ee_{8(-24)}$ to its subalgebra $\ee_{6(-26)}$ by considering
the following generators.
\begin{itemize}
\item
The 28 elements $\eA_p$, $\eG_p$, $\eS_p$, $\eD_p$, with $p\in\UU$, generate
$\so(8)$.
\item
After adding the 24 elements $\eX_q$, $\eY_q$, $\eZ_q$, with
$q\in\UU\cup\{1\}$, the resulting 52 elements generate~$\ff_4$.
\item
The 24 elements $\eX_{qL}$, $\eY_{qL}$, $\eZ_{qL}$ are boosts, with
$q\in\UU\cup\{1\}$ (and $L\in\CC'\subset\OO'$ with $L^2=+1$).
\item
The two remaining elements, both boosts, are~$\eD_L$ and~$\eS_L$.
\end{itemize}
Since we have only used indices in $\CC'\otimes\OO$, it is straightforward to
verify that these $28+24+26=78$ elements close, and have the correct
commutators; we have generated $\ee_{6(-26)}\isom\su(3,\CC'\otimes\OO)$.  With
this labeling, it is easy to distinguish boosts from rotations; boosts contain
$L$ in their indices, while rotations do not.

The centralizer of $\ee_6$ in $\ee_8$ is $\su(3)$, in this case the real form
of $\su(3)$ in $\gg_{2(2)}=\der(\OO')$ that fixes $L$.  That real form is
$\sl(3,\RR)$, given explicitly by
\begin{equation}
\sl(3,\RR)
  = \langle \eA_L,\eG_L,\eA_I,\eA_J,\eA_K,\eA_{KL},\eA_{JL},\eA_{IL} \rangle
\end{equation}
where the angled brackets denote the span of the given generators.

Which generators of $\ee_8$ are left?  That's easy: all elements containing a
split imaginary unit other than $L$.  So for each $Q\in\{I,J,K,KL,JL,IL\}$, we
have the 27 elements $\eX_{pQ}$, $\eY_{pQ}$, $\eZ_{pQ}$, $\eD_Q$, $\eS_Q$, and
$\eG_Q$, for a total of $6\times27=162$ elements left over, thus accounting
for all $78+8+162=248$ elements of $\ee_8$.
%

Since the only 27-dimensional representation of $\ee_6$ is the Albert algebra,
it is clear from representation theory that the induced action (using
commutators in $\ee_8$) of $\ee_6$ on each of the $27$s in the decomposition
must reproduce the action of $\ee_6$ on the Albert algebra.  Since this action
can not change the value of $Q$ in the elements given above, we conclude that
we have six copies of the Albert algebra, labeled by $Q$.

However, we expect these $27$s to be simultaneous representations of both
$\ee_6$ and $\sl(3,\RR)$, which constrains the choice of labels.  We therefore
digress briefly to discuss the structure of $\sl(3,\RR)$.

\subsection{\boldmath The structure of $\sl(3,\RR)$}
\label{sl3R}

The minimal representation of $\gg_{2(2)}$ is 7-dimensional---the imaginary
split octonions.  The subalgebra $\sl(3,\RR)\subset\gg_{2(2)}=\der(\OO')$
corresponds at the group level to fixing $L$, so $\sl(3,\RR)$ acts on
$\CC'_\perp\isom\Im\HH\oplus(\Im\HH)L$, a three-dimensional vector space over
$\CC'$.  An explicit $3\times3$ (over $\CC'$) matrix representation of
$\sl(3,\RR)$ is given in Appendix~\ref{sl3Rmat}, but the matrices given there
are \textit{not} elements of $\ee_8$, as they ignore the imaginary elements
of~$\CC'_\perp$.

Nonetheless, choosing $\eA_L$ and $\eG_L$ to generate the Cartan subalgebra of
$\gg_{2(2)}$, it is straightforward to determine their eigenvectors in $\OO'$,
which are $\{I\pm IL,J\pm JL, K\pm KL\}$,
%
%
suggesting that we use these six null elements of $\OO'$ to label the Jordan
algebras.

\subsection{\boldmath The action of $\ee_6$ on an Albert algebra}
\label{action}

We return to the six 27-dimensional representations introduced in
Section~\ref{e6gen}, whose generators have the form
\begin{equation}
\HQ = \langle \eX_{pQ},\eY_{pQ},\eZ_{pQ},\eD_Q,\eS_Q,\eG_Q \rangle
\end{equation}
and show that each such representation is equivalent to the usual Hermitian
representation of the Albert algebra.  However, we now use the eigenvectors of
$\sl(3,\RR)$ as labels, thus taking $2Q\in\{I\pm IL,J\pm JL,K\pm KL\}$.

In this signature, $\ee_6$ contains 52 rotations, making up $\ff_4$, and 26
boosts.  We begin with the action of $\ff_4$ on $\HQ$, then consider the
boosts, paying particular attention to the ``trace'' element, $\eG_Q$.  Since
$\ff_4$ is the centralizer of $\gg_{2(2)}$ in $\ee_8$, $\ff_4$ must commute
with $\eG_Q\in\gg_{2(2)}$.  The remaining 26 elements in $\HQ$ are
anti-Hermitian matrices having $Q$ as an overall factor; factoring out the $Q$
yields a Hermitian matrix, and in fact a tracefree element of the Albert
algebra $\HHH$---precisely the minimal representation of $\ff_4$.  We must
therefore show that the $\ee_8$ action reproduces the expected action
of~$\ff_4$ on these matrices.

A basis element of $\ff_4$ is either an octonionic anti-Hermitian tracefree
matrix, or can be obtained as the commutator of two such matrices; it
therefore suffices to consider the action of such matrices.  So let
$\AAA\in\sa(3,\OO)\subset\ff_4$, which acts on a tracefree element
$\PPP\in\HHH$ via
\begin{equation}
\PPP \longmapsto \AAA\PPP + \PPP\AAA^\dagger = [\AAA,\PPP]
\end{equation}
since $\AAA$ is anti-Hermitian.  Since $Q$ commutes with these commutators,
that is,
\begin{equation}
[\AAA,Q\PPP] = Q[\AAA,\PPP]
\end{equation}
we have shown that the $\ee_8$ action reduces to the $\ff_4$ action on the
tracefree part of $\HQ$.  But we showed above that $\ff_4$ commutes with
$\eG_Q$, as it should if $\eG_Q$ corresponds to a trace (that is, to a
multiple of the identity element in the Albert algebra).  The $\ee_8$ action
therefore does reproduce the action of $\ff_4$ on $\HQ$, as expected.

What about the action of boosts, that is, of the remaining elements of
$\ee_6\subset\ee_8$ on $\HQ$?

In traditional language, boosts are generated by the 26 \textit{Hermitian}
tracefree $3\times3$ octonionic matrices.  Let $\BBB$ be one of these
matrices.  Then $\BBB L\in\ee_8$, and is anti-Hermitian.  The $\ee_6$ action
of $\BBB$ on a tracefree element $\PPP\in\HHH$ is given by
\begin{equation}
\PPP \longmapsto \BBB\PPP + \PPP\BBB^\dagger = \{\BBB,\PPP\}
\end{equation}
whereas the $\ee_8$ action of $\BBB L$ on $Q\PPP$ is
\begin{equation}
[\BBB L,Q\PPP]
  = (\BBB L)(Q \PPP) - (Q\PPP)(\BBB L) = (LQ) \{\BBB,\PPP\}
\label{e8boost}
\end{equation}
since $L$ anticommutes with $Q$ (and $\OO'$ commutes with $\OO$).
Furthermore, $(1\pm L)/2$ are projection operators in $\OO'$, with $LQ=\pm Q$.

We're almost there, as we have shown that the $\ee_8$ action of the boost
$\BBB L$ reproduces the $\ee_6$ action of the boost $\BBB$ up to sign.  Since
the relative sign between the action of boosts and rotations distinguishes
between the $27$ representation of $\ee_6$ and its dual, the $\bar{27}$, we
have shown that the six copies of $\HQ$ consist of three copies of each.

There is however one remaining subtlety, as we have only considered the action
of boosts on tracefree elements of the Albert algebra.  So we must now check
the action of boosts when the nested element $\eG_Q$ is involved.
Since~\eqref{Gp} holds unchanged over $\OO'$, we have
\begin{equation}
\eG_Q = -2 (QL\circ L)\,\III = \mp2 (Q\circ L)\,\III
\end{equation}
(since $QL=\pm Q$).  We can therefore use~\eqref{Dc}--\eqref{Fc} to express
$\eD_{Q,L}$, $\eE_{Q,L}$, and $\eF_{Q,L}$ in terms of the single nested
element $\eG_Q$ together with tracefree elements of $\HQ$.  Since the only
nested factor that can occur in~\eqref{e8boost} is $Q\circ L$, and since this
factor occurs only in $\eD_{Q,L}$, $\eE_{Q,L}$, $\eF_{Q,L}$ and elements built
from them,~\eqref{e8boost} does indeed show that the $\ee_8$ action of $\BBB
L$ reduces to the $\ee_6$ action of $\BBB$ when acting on \textit{traceless}
elements $Q\PPP$ and $\PPP$, respectively.

The final step is to consider
\begin{align}
[\BBB L,\frac12\eG_Q]
  &= [\BBB L,\mp(Q\circ L)\,\III] \nonumber\\
  &= \bigl( \mp L \circ (Q\circ L) \pm (Q\circ L)\circ L \bigr) \BBB
	\nonumber\\
  &= (\pm Q \pm Q ) \BBB
   = \pm Q \,\{\BBB,\III\}
   = QL \,\{\BBB,\III\}
\end{align}
which shows that $\BBB L$ acting on $\frac12\eG_Q$ in $\ee_8$ is the same (up
to sign) as $Q$ times the action of $\BBB$ on $\III$ in $\ee_6$, with the sign
again indicating which of the two representations of $\ee_6$ we have, which are
dual to each other.

We have shown that the action of $\ee_6\subset\ee_8$ on $\HQ$ is equivalent to
an action of $\ee_6$ on $\HHH$.  By providing a matrix language for discussing
nested elements such as $\eG_Q$, we are able to describe the $27$s of $\ee_6$
as generalized matrices in $\ee_8$.  We have further shown that these
anti-Hermitian matrices are just the usual representation of the Albert
algebra in terms of $3\times3$ Hermitian matrices, multiplied by an overall
``label'' $Q\in\OO'$, where ``$Q\III$'' is interpreted as $-\frac12\eG_Q$.
This interpretation is a key piece of our generalized matrix formalism.

\setlength{\textheight}{9.6in}
\section{Other Cases}
\label{discussion}

The construction of $\ee_8$ given in Section~\ref{e8} works for all three real
forms of $\ee_8$, namely
\begin{align}
\ee_{8(-248)} &\isom \su(3,\OO\otimes\OO) ,\\
\ee_{8(-24)} &\isom \su(3,\OO'\otimes\OO) ,\\
\ee_{8(8)} &\isom \su(3,\OO'\otimes\OO') ,
\end{align}
depending on the choice of division algebras used.
%
%
%
%
In this paper, we have emphasized the split real forms $\ee_{6(-26)}\subset\ee_{8(-24)}$.
However, there are in fact three real forms of $\ee_6$ contained in
$\ee_{8(-24)}\isom\su(3,\OO'\otimes\OO)$, namely
$\ee_{6(-26)}\isom\su(3,\CC'\otimes\OO)$,
$\ee_{6(2)}\isom\su(3,\OO'\otimes\CC)$, and compact
$\ee_{6(-78)}\isom\su(3,\CC\otimes\OO)$; the split real form
$\ee_{6(6)}\isom\su(3,\CC'\otimes\OO')$ is contained in the split real form
$\ee_{8(8)}\isom\su(3,\OO'\otimes\OO')$.  (The remaining real form of $\ee_6$,
namely $\ee_{6(-14)}\isom\sl(2,1,\OO)$, is a subalgebra of
$\ee_{8(-24)}\isom\su(3,\OO'\otimes\OO)$, but is not itself of the form
$\su(3,\KK_1\otimes\KK_2)$.)

The representation theory is slightly different in the compact case, where we
expect the decomposition of $\ee_8$ over $\ee_6\oplus\su(3)$ to be
\begin{equation}
248 = 78 + 3\times54 + 8
\end{equation}
with three \textit{complex} Albert algebras, rather than six \textit{real}
Albert algebras.  This behavior is easy to track in our notation, since the
separation into six real Albert algebras relied on the existence of projection
operators $1\pm L$, which only exist in $\OO'$.

However, the real form $\ee_{8(-24)}\isom\su(3,\OO'\otimes\OO)$ appears to be
the most interesting from the point of view of particle physics.  Manogue and
Dray~\cite{Denver,York,Spin} previously proposed using $\HHH$ to
represent three generations of leptons; a quite different interpretation
leading to three generations was given in~\cite{Octions}.  However, in both
cases, the spinor elements of the six copies of $\HQ$ do indeed correspond to
the three colors and three anticolors of the Standard Model.

\section{Discussion}

In summary, although the decomposition of (real forms of) $\ee_8$ into (real
forms of) $\ee_6$ and $\su(3)$ and their representations is well known, it is
rarely emphasized that the adjoint action of $\ee_8$ \textit{must} therefore
reproduce the usual action of $\ee_6$ on its minimal representations, namely
the Albert algebra.  We have verified this expectation explicitly, using a
newly-discovered matrix formalism for $\ee_8$.  This (re)interpretation of
known actions on given representations in terms of the adjoint action of a
larger symmetry group lies at the heart of our recent proposal~\cite{Octions}
to embed \textit{all} of the objects of the Standard Model within~$\ee_8$.

Along the way, we have introduced a new representation of the Albert algebra
as \textit{anti-Hermitian} matrices in $\ee_8$ that combine the traditional
Hermitian matrices over $\OO$ with imaginary labels in $\OO'$.  Not only do
matrix commutators in $\ee_8$ reproduce the action of $\ee_6$ on its
minimal representation, as shown here, but they also reproduce the structure
of the Albert algebra as a Jordan algebra, as discussed in~\cite{MagicE7}.
This representation may therefore be useful in other recent efforts to use the
Albert algebra to describe physical phenomena, such as~\cite{Singh22}, or,
more generally, in other approaches involving a distinguished copy of $\ee_6$
in $\ee_8$, such as~\cite{Chirality22}.

\acknowledgments{
This work was supported in part by the John Templeton Foundation under grant
number 34808, by FQXi, and by the Institute for Advanced Study.
}

\appendix

\section{\boldmath A preferred basis for $\so(8)$}
\label{so8}

\begin{table}
\begin{center}
\begin{tabular}{|c|c|c|c|}
\hline
$p$ & $a(p)$ & $b(p)$ & $c(p)$ \\
\hline
\hline
$i$  & $j,k $ & $k\ell,j\ell$ & $\ell,i\ell$ \\
$j$  & $k,i $ & $i\ell,k\ell$ & $\ell,j\ell$ \\
$k$  & $i,j $ & $j\ell,i\ell$ & $\ell,k\ell$ \\
$k\ell$ & $j\ell,i$ & $j,i\ell $ & $k,\ell $ \\
$j\ell$ & $i\ell,k $ & $i,k\ell$ & $j,\ell $ \\
$i\ell$ & $k\ell,j$ & $k,j\ell $ & $i,\ell $ \\
$\ell$  & $i\ell,i$ & $j\ell,j $ & $k\ell,k$ \\
\hline
\end{tabular}
\caption{Each pair $a(p)$, $b(p)$, $c(p)$ generates a quaternionic subalgebra
of $\OO$ that contains $p$.}
\label{AGS}
\end{center}
\end{table}

We have constructed three different bases of $\so(8)$ in Section~\ref{e8}.
Here, we construct yet another basis of $\so(8)$, which is adapted to the
preferred $\gg_2$ subalgebra of $\so(8)$ that preserves the octonionic
multiplication table.

Each of these bases contains seven ordinary matrices (of the form $\eD_p$,
$\eE_p$, or $\eF_p$), and $\binom72=21$ ``nested'' matrices~(of the form
$\eD_{p,q}$, $\eE_{p,q}$, or $\eF_{p,q}$), constructed using~$\circ$, and
where $p\ne q\in\UU$.  As discussed in~\cite{Structure}, for each $p\in\UU$
there are three ordered pairs of elements of $\UU$ that ``point'' to $p$ in
the Fano plane (Figure~\ref{omult3}), corresponding to 3~quaternionic
subalgebras containing $p$; see Table~\ref{AGS} (which differs slightly from
that given in~\cite{Structure}).  In this table, $a(p)$, $b(p)$, $c(p)$
represent the three independent choices of pairs of imaginary units in $\OO$
that generate a quaternionic subalgebra of $\OO$ that contains $p$.  We use
these pairs to construct a new basis for the nested elements $\eD_{p,q}$,
namely
\begin{align}
-\eA_p &= \eD_{a(p)} - \eD_{b(p)} \label{Adef} ,\\
-\eG_p &= \eD_{a(p)} + \eD_{b(p)} -2 \eD_{c(p)} \label{Gdef} ,\\
-\eS_p &= \eD_{a(p)} + \eD_{b(p)} + \eD_{c(p)} \label{Sdef}
\end{align}
with a mild but hopefully obvious abuse of notation.

This basis has some useful properties, which can be verified by direct
computation.  First of all, the 14 elements $\{\eA_p,\eG_p\}$ generate
$\gg_2$.  Furthermore, the expressions given above for $\eA_p$ and $\eG_p$ are
\textit{type-independent}, in the sense that $\eD$ can be replaced by $\eE$ or
$\eF$ in~\eqref{Adef} and~\eqref{Gdef} without changing the result, a property
that we refer to as \textit{strong triality}.  Triality also turns out to
imply that
\begin{equation}
\eS_p
  = \eE_p - \eF_p
  = \begin{pmatrix}p& 0& 0\\ 0& p& 0\\ 0& 0& -2p\\\end{pmatrix}
\end{equation}
which shows how to construct $\eS_p$ without nesting.  Finally, since
\begin{equation}
\eD_p + \eE_p + \eF_p = 0
\end{equation}
we can express any of the single- or double-subscript $\eE$s and $\eF$s in terms
of $\{\eA_p,\eG_p,\eS_p,\eD_p\}$, which is the promised basis of $\so(8)$.
The $\eA$s and $\eG$s are the only nested elements, and generate $\gg_2$; adding
the $\eS$s generates $\so(7)$; adding the $\eD$s yields all of $\so(8)$.
Explicitly,
\begin{align}
\eE_p &= -\frac12\, (\eD_p-\eS_p) ,\\
\eF_p &= -\frac12\, (\eD_p+\eS_p) ,
\end{align}
and for instance
\begin{align}
\eD_{c(p)}
 &= -\frac13\, (\eE_p-\eF_p-\eG_p)
  = \frac13\, (\eG_p-\eS_p) ,\label{Dc}\\
\eE_{c(p)}
 &= -\frac13\, (\eF_p-\eD_p-\eG_p)
  = \frac16\, (\eS_p+2\eG_p+3\eD_p) ,\label{Ec}
\\
\eF_{c(p)}
 &= -\frac13\, (\eD_p-\eE_p-\eG_p)
  = \frac16\, (\eS_p+2\eG_p-3\eD_p) ,\label{Fc}
\end{align}
with similar expressions holding for the remaining nested elements.

We emphasize that our new basis of $\so(8)$ is of type~I; commutators with
elements of types~II and~III require rewriting our basis in terms of elements
of the appropriate type.  As noted above, this step is trivial for the $\eA$s
and $\eG$s, but a bit of algebra is required in the remaining cases.

We note for future reference that~\eqref{Dc}--\eqref{Fc} imply that
\begin{equation}
\eD_{c(p)} + \eE_{c(p)} + \eF_{c(p)} = \eG_p ;
\label{DEF}
\end{equation}
a similar construction yields
\begin{equation}
(\eD_{a(p)} + \eE_{a(p)} + \eF_{a(p)}) - (\eD_{b(p)} + \eE_{b(p)} + \eF_{b(p)})
 = -3\eA_p .
\end{equation}
Comparing~\eqref{DEF} with~\eqref{comp}, we see that each $\eG_p$ can be
written as a multiple of the identity matrix~$\III$, namely
\begin{equation}
\eG_p = -2(p\ell\circ\ell)\,\III
\label{Gp}
\end{equation}
for $\ell\ne p\in\UU$, and that
\begin{equation}
\eG_\ell = 2(k\ell\circ k)\,\III
\end{equation}
with similar expressions holding for $\eA_p$.

\section{\boldmath A matrix representation of $\sl(3,\RR)$}
\label{sl3Rmat}

The Lie algebra $\sl(3,\RR)$ is a real form of $\su(3)$, with five boosts and
three rotations.  Its generators are usually written as $3\times3$ real
matrices, in which case the rotations are anti-Hermitian and the boosts are
Hermitian.  Here, we interpret $\sl(3,\RR)$ instead as $\su(3,\CC')$, which
has the effect of multiplying the boosts by $L\in\CC'$, resulting in all eight
generators being anti-Hermitian.  An explicit basis is given by
\begin{align}
\frac12\,\eA_I
  &= \begin{pmatrix} 0& 0& 0\\ 0& 0& -1\\ 0& 1& 0\\\end{pmatrix} ,
\qquad
\frac12\,\eA_{IL}
  = \begin{pmatrix} 0& 0& 0\\ 0& 0& L\\ 0& L& 0\\\end{pmatrix} ,
\\
\frac12\,\eA_J
  &= \begin{pmatrix} 0& 0& 1\\ 0& 0& 0\\ -1& 0& 0\\\end{pmatrix} ,
\qquad
\frac12\,\eA_{JL}
  = \begin{pmatrix} 0& 0& L\\ 0& 0& 0\\ L& 0& 0\\\end{pmatrix} ,
\\
\frac12\,\eA_K
  &= \begin{pmatrix} 0& -1& 0\\ 1& 0& 0\\ 0& 0& 0\\\end{pmatrix} ,
\qquad
\frac12\,\eA_{KL}
  = \begin{pmatrix} 0& L& 0\\ L& 0& 0\\ 0& 0& 0\\\end{pmatrix} ,
\\ 
\frac12\,\eA_L
  &= \begin{pmatrix} L& 0& 0\\ 0& -L& 0\\ 0& 0& 0\\\end{pmatrix} ,
\qquad
\frac12\,\eG_L
  = \begin{pmatrix} L& 0& 0\\ 0& L& 0\\ 0& 0& -2L\\\end{pmatrix} .
\end{align}
The structure constants for this Lie algebra are real; it is indeed a copy
of~$\sl(3,\RR)$.

We emphasize that these matrices are the the analogs of the Gell-Mann matrices
for $\su(3)$, but are \textit{not} matrix representations of $\sl(3,\RR)$ in
$\ee_8$.
To reinterpret $\sl(3,\RR)$ as a subalgebra of $\ee_{8(-24)}$ requires
reversing the identification given in Section~\ref{sl3R}, thus identifying
the three-dimensional (over $\CC'$) vector space on which these matrices act
with the six-dimensional (over $\RR$) vector space $\CC'_\perp\subset\OO'$.

\bibliographystyle{unsrt}
\bibliography{octo,octo2}

\end{document}